\documentclass[12pt]{amsart}
\usepackage[cp1251]{inputenc}
\usepackage[english,russian]{babel}
\usepackage{amsmath}
\usepackage{amssymb}
\usepackage{amsfonts}
\title{Potentials for Elliptic Boundary Value Problems in Cones}

\author{Vladimir Vasilyev }

\date{}

\begin{document}
\renewcommand{\refname}{References}
\maketitle

{\bf Abstract.} We consider an elliptic pseudo differential equation in a multi-dimensional cone and starting wave factorization concept we add some boundary conditions.
For the simplest cases explicit formulas for solution are given like layer potentials for classical case.



\section{Introduction}

How potentials are constructed for boundary value problems? One takes a fundamental solution of corresponding differential operator {\it in whole space ${\bf R}^m$}, and with its help constructs the potentials
according to boundary conditions. Further, one studies their boundary properties, and with the help of potentials reduces the boundary value problem to an equivalent integral equation on the boundary. The formulas
for integral representation of solution of the boundary value problem one has for separate cases only (ball, half-space, such places, where one has explicit form for Green function). Thus, an ideal result for boundary value problem
even with smooth boundary is its reduction to equivalent Fredholm equation and obtaining the existence and uniqueness theorem (without knowing how the solution  looks) \cite{A,F,K,MMT,HW}. I would like to show, that potentials can arise from another
point of view, without using fundamental solution, but using factorization idea and they obviously must take into account the boundary geometry. Smooth boundary is locally a hyper-plane (there is Poisson formula for the
Dirichlet problem, see also \cite{E}),
first type of non-smooth boundary is conical surface.\footnote{
This work was completed when the author was a DAAD stipendiat and hosted in Institute of Analysis and Algebra, Technical University of Braunschweig.}

\section{Preliminaries}

Let's go to studying solvability of pseudo differential equations \cite{V1,V2,V3}
\begin{equation}\label{1}
(Au_+)(x)=f(x),~x\in C^a_+,
\end{equation}
in the space $H^s(C^a_+)$, where $C^a_+$ is $m$-dimensional cone $C^a_+=\{x\in{\bf R}^m: x=(x_1,...,x_{m-1},x_m), x_m>a|x'|, a>0\},~x'=(x_1,...,x_{m-1})$, $A$ is pseudo differential operator ($\tilde u$ denotes the Fourier transform of $u$)
\[
u(x)\longmapsto\int\limits_{{\bf R}^m}e^{ix\cdot\xi}A(\xi)\tilde u(\xi)d\xi,~~x\in{\bf R}^m,
\] 
 with the symbol$A(\xi)$ satisfying the condition
$$
c_1\leq|A(\xi)(1+|\xi|)^{-\alpha}|\leq c_2.
$$
(Such symbols are elliptic \cite{E} and have the order $\alpha\in{\bf R}$ at infinity.)

By definition, the space $H^s(C^a_+)$ consists of distributions from $H^s({\bf R}^m)$, which support belongs to $\overline{C^a_+}$. The norm in the space $H^s(C^a_+)$ is induced by the norm from $H^s({\bf R}^m)$. The right-hand side $f$ is
chosen from the space $H^{s-\alpha}_0(C^a_+)$, which is space of distributions
 $S'(C^a_+)$, admitting the continuation on $H^{s-\alpha}({\bf R}^m)$. The norm in the space  $H^{s-\alpha}_0(C^a_+)$ is defined
$$
||f||^+_{s-\alpha}=\inf ||lf||_{s-\alpha},
$$
where {\it infimum} is chosen from all continuations $l$.

Further, let's define a special multi-dimensional singular integral by the formula
\[
(G_mu)(x)=\lim\limits_{\tau\to 0+}\int\limits_{\bf{R}^m}
\frac{u(y',y_m)dy'dy_m}{\left(|x'-y'|^2-a^2(x_m-y_m+i\tau)^2\right)^{m/2}}
\]
(we omit a certain constant, see \cite{V1}). Let us recall, this operator is multi-dimensional analogue of one-dimensional Cauchy type integral, or Hilbert transform.

We  also need some notations before definition.

The symbol $\stackrel{*} {C^a_+}$ denotes a conjugate cone for $C^a_+$:
\[
\stackrel{*} {C^a_+}=\{x\in{\bf R}^m: x=(x',x_m), ax_m>|x'|\},
\]
$C^a_-\equiv -C^a_+,~T(C^a_+)$ denotes radial tube domain over the cone $C^a_+$, i.e. domain in a complex space ${\bf C}^m$ of type ${\bf R}^m+iC^a_+$.

To describe the solvability picture for the equation \eqref{1} we will introduce the following

{\bf Definition.}
{\it Wave factorization for the symbol $A(\xi)$ is called its representation in the form
$$
A(\xi)=A_{\neq}(\xi)A_=(\xi),
$$
where the factors $A_{\neq}(\xi),A_=(\xi)$ must satisfy the following conditions:

1) $A_{\neq}(\xi),A_=(\xi)$ are defined for all admissible values $\xi\in{\bf R}^m$, without may be, the points $\{\xi\in{\bf R}^m:|\xi'|^2=a^2\xi^2_m\}$;

2) $A_{\neq}(\xi),A_=(\xi)$ admit an analytical continuation into radial tube domains\\ $T(\stackrel{*} {C^a_+}),T(\stackrel{*} {C^a_-})$ respectively with estimates
$$
|A_{\neq}^{\pm 1}(\xi+i\tau)|\leq c_1(1+|\xi|+|\tau|)^{\pm\ae},
$$
$$
|A_{=}^{\pm 1}(\xi-i\tau)|\leq c_2(1+|\xi|+|\tau|)^{\pm(\alpha-\ae)},~\forall\tau\in\stackrel{*} {C^a_+}.
$$

The number $\ae\in{\bf R}$ is called index of wave factorization.}

The class of elliptic symbols admitting the wave factorization is very large. There are the special chapter in the book \cite{V1} and the paper \cite{V4} devoted to this question, there are examples also for certain operators of mathematical physics.

Everywhere below we will suppose that the  mentioned wave factorization does exist, and the sign $\sim$ will denote the Fourier transform, particularly $\tilde H(D)$ denotes the Fourier image of the space $H(D)$.

\section{Scheme in details}\label{sec.3}

Now we will consider the equation \eqref{1}
for the case $\ae-s=n+\delta, n\in{\bf N}, |\delta|<1/2$, only. A general solution can be constructed by the following way. We choose an arbitrary continuation $lf$ of the right-hand side on a whole space $H^{s-\alpha}({\bf R}^m)$ and introduce
\[
u_-(x)=(lf)(x)-(Au_+)(x).
\]

After wave factorization for the symbol $A(\xi)$ with preliminary Fourier transform we write
$$
A_{\neq}(\xi)\tilde u_+(\xi)+A^{-1}_=(\xi)\tilde u_-(\xi)=A^{-1}_=(\xi)\tilde{lf}(\xi).
$$

One can see that $A^{-1}_=(\xi)\tilde{lf}(\xi)$  belongs to the space $\tilde H^{s-\ae}({\bf R}^m)$, and if we choose the polynomial $Q(\xi)$, satisfying the condition
$$
|Q(\xi)|\sim(1+|\xi|)^n,
$$
then $Q^{-1}(\xi)A^{-1}_=(\xi)\tilde{lf}(\xi)$ will belong to the space $\tilde H^{-\delta}({\bf  R}^m)$.

Further, according to the theory of multi-dimensional Riemann problem \cite{V1},
we can decompose the last function on two summands (jump problem):
$$
Q^{-1}A^{-1}_=\tilde{lf}=f_++f_-,
$$
where $f_+\in\tilde H(C^a_+), f_-\in\tilde H({\bf R}^m\setminus C^a_+).$

So, we have
$$
Q^{-1}A_{\neq}\tilde u_++Q^{-1}A_{=}^{-1}\tilde u_-=f_++f_-,
$$
or
$$
Q^{-1}A_{\neq}\tilde u_+-f_+=f_--Q^{-1}A_{=}^{-1}\tilde u_-
$$

In other words,
$$
A_{\neq}\tilde u_+-Qf_+=Qf_--A_=^{-1}\tilde u_-.
$$

The left-hand side of the equality belongs to the space $\tilde H^{-n-\delta}(C^a_+)$, and right-hand side is from $\tilde H^{-n-\delta}({\bf R}^m\setminus C^a_+)$, hence
$$
F^{-1}(A_{\neq}\tilde u_+-Qf_+)=F^{-1}(Qf_--A_=^{-1}\tilde  u_-),
$$
where the left-hand side belongs to the space $ H^{-n-\delta}(C^a_+)$, and the right-hand side belongs to the space $ H^{-n-\delta}({\bf R}^m\setminus C^a_+)$,
that's why we conclude immediately that it is distribution supported on $\partial C^a_+$.

the main tool now is to define the form of the distribution.

Let's denote $T_a$ the bijection operator transferring  $\partial C^a_+$  into hyperplane $x_m=0$, more precisely, it is transformation
${\bf R}^m\longrightarrow{\bf R}^m$ of the following type
$$
\left\{
\begin{array}{cccc}
t_1=x_1,\\
............\\
t_{m-1}=x_{m-1},\\
t_m=x_m-a|x'|.
\end{array}
\right.
$$

Then the function
$$
T_aF^{-1}(A_{\neq}\tilde u_+-Qf_+)
$$
will be supported on the hyperplane $t_m=0$ and belongs to $ H^{-n-\delta}({\bf R}^m)$. Such distribution is a linear span of Dirac mass-function and its derivatives \cite{GS} and looks as the following sum
$$
\sum\limits_{k=0}^{n-1}c_k(t')\delta^{(k)}(t_m).
$$

It is left to think, what is operator $T_a$ in Fourier image. Explicit calculations give simple answer:
$$
FT_au=V_a\tilde u,
$$
where $V_a$ is something like a pseudo differential operator with symbol $e^{-ia|\xi'|\xi_m}$, and further one can construct the general solution for our pseudo differential equation \eqref{1}.

We need some connections between the Fourier transform and the operator $T_a$:
$$
(FT_au)(\xi)=\int\limits_{{\bf R}^m}e^{-ix\cdot\xi}u(x_1,...,x_{m-1},x_m-a|x'|)dx=
$$
$$
=\int\limits_{{\bf R}^m}e^{-iy'\xi'}e^{-i(y_m+a|y'|)\xi_m}u(y_1,...,y_{m-1},y_m)dy=
$$
$$
=\int\limits_{{\bf R}^{m-1}}e^{-ia|y'|\xi_m}e^{-iy'\xi'}\hat{u}(y_1,...,y_{m-1},\xi_m)dy',
$$
where $\hat{u}$ denotes the Fourier transform on the last variable, and Jacobian is

$$
\frac{D(x_1,x_2,...,x_m)}{D(y_1,y_2,...,y_m)}=\left |
\begin{array}{ccccc}
1&0&\cdots0&0&\\
0&1&\cdots0&0&\\
0&0&\cdots1&0\\
&\cdots\cdots\cdots\cdots\cdots\\
\frac{ay_1}{|y'|}&\frac{ay_2}{|y'|}&....\frac{ay_{m-1}}{|y'|}&1

\end{array}
\right |=1.
$$

If we define a pseudo differential operator by the formula
\[
(Au)(x)=\int\limits_{{\bf R}^m}e^{ix\xi}A(\xi)\tilde u(\xi)d\xi,
\]
and the direct Fourier transform
$$
\tilde u(\xi)=\int\limits_{{\bf R}^m}e^{-ix\xi}u(x)dx,
$$
then we have the following relation formally at least
\begin{equation}\label{2}
(FT_au)(\xi)=\int\limits_{{\bf R}^{m-1}}e^{-ia|y'|\xi_m}e^{-iy'\xi'}\hat{u}(y_1,...,y_{m-1},\xi_m)dy.
\end{equation}

In other words, if we denote the $(m-1)$-dimensional Fourier transform ($y'\to\xi'$ in distribution sense) of function $e^{-ia|y'|\xi_m}$ by $E_a(\xi',\xi_m)$, then
the formula \eqref{2} will be the following
\[
(FT_au)(\xi)=(E_a*\tilde u)(\xi),
\]
where the sign $*$ denotes a convolution for the first $m-1$ variables, and the multiplier for the last variable $\xi_m$. Thus, $V_a$ is a combination of a convolution operator and the multiplier with the kernel $E_a(\xi',\xi_m)$. It is very simple operator, and it is bounded in Sobolev-Slobodetskii spaces $H^s({\bf R}^m)$.

Notice that distributions supported on conical surface and their Fourier transforms were considered in \cite{GS}, but the author didn't find the multi-dimensional analogue of theorem on a distribution supported in a
single point in all issues of this book.

{\bf Remark 1.}
{\it  One can wonder why we can't use this transform in the beginning to reduce the conical situation \eqref{1} to hyper-plane one, and then to apply Eskin's technique \cite{E}. Unfortunately, it's impossible, because
$T_a$ is non-smooth transformation, but even for smooth transformation we obtain the same operator $A$ with some additional compact operator. Obtaining the invertibility conditions for such operator is a very serious problem.}

\section{General solution}\label{sec3}

The following result is valid (it follows from considerations of Sec.3).

{\bf Theorem.}\label{t1}
A general solution of the equation \eqref{1} in Fourier image is given by the formula
\[
\tilde u_+(\xi)=A^{-1}_{\neq}(\xi)Q(\xi)G_mQ^{-1}(\xi)A^{-1}_=(\xi)\tilde{lf}(\xi)+
\]
\[
+A^{-1}_{\neq}(\xi)V_{-a}F\left(\sum\limits_{k=1}^nc_k(x')\delta^{(k-1)}(x_m)\right),
\]
where $c_k(x')\in H^{s_k}({\bf R}^{m-1})$ are arbitrary functions, $s_k=s-\ae+k-1/2,~k=1,2,...,n,$ $lf$ is an arbitrary continuation $f$ on $H^{s-\alpha}({\bf R}^m)$.

Starting this representation one can suggest different statements of boundary value problems for the equation \eqref{1}.

\section{Boundary conditions: simplest variant, the Dirichlet condition}\label{sec4}

Let's consider a very simple case, when $f\equiv 0,~a=1,~n=1$. Then the formula from theorem  takes the form
$$
\tilde u_+(\xi)=A^{-1}_{\neq}(\xi)V_{-1}\tilde c_0(\xi').
$$

We consider the following construction separately. According to the Fourier transform our solution is
\[
u_+(x)=F^{-1}\{A^{-1}_{\neq}(\xi)V_{-1}\tilde c_0(\xi')\}.
\]

Let's suppose we choose the Dirichlet boundary condition on $\partial C^1_+$ for unique identification of an unknown function $c_0$, i.e.
$$
(Pu)(y)=g(y),
$$
where $g$ is given function on $\partial C^1_+$, $P$ is restriction operator on the boundary,
so we know the solution on the boundary $\partial C^1_+$.

Thus,
\[
T_1u(x)=T_1F^{-1}\{A^{-1}_{\neq}(\xi)V_{-1}\tilde c_0(\xi')\},
\]
so we have
\begin{equation}\label{3}
FT_1u(x)=FT_1F^{-1}\{A^{-1}_{\neq}(\xi)V_{-1}\tilde c_0(\xi')\}=V_1\{A^{-1}_{\neq}(\xi)V_{-1}\tilde c_0(\xi')\},
\end{equation}
and we know $(P'T_1u)(x')\equiv v(x')$,
where $P'$ is the restriction operator on the hyperplane $x_m=0$.

The relation between the operators $P'$ and $F$ is well-known \cite{E}:
$$
(FP'u)(\xi')=\int\limits_{-\infty}^{+\infty}\tilde u(\xi',\xi_m)d\xi_m.
$$

Returning to the formula \eqref{3} we obtain the following
\begin{equation}\label{4}
\tilde v(\xi')=\int\limits_{-\infty}^{+\infty}\{V_1\{A^{-1}_{\neq}(\xi)V_{-1}\tilde c_0(\xi')\}\}(\xi',\xi_m)d\xi_m,
\end{equation}
where $\tilde v(\xi')$ is given function. Hence, the equation \eqref{4} is an integral equation for determining $c_0(x')$.

The Neumann boundary condition leads to analogous integral equation (see below).

\section{Conical potentials}\label{sec5}

Let's consider the particular case: $f\equiv 0, n=1$. The formula for general solution of the equation \eqref{1} takes the form
$$
\tilde u_+(\xi))=A_{\neq}^{-1}(\xi)V_{-a}F\{c_0(x')\delta^{(0)}(x_m)\},
$$
and further after Fourier transform (for simplicity we write $\tilde c$ instead of $V_{-1}\tilde c_0$)
\begin{equation}\label{5}
\tilde u_+(\xi)=A_{\neq}^{-1}(\xi)\tilde c(\xi'),
\end{equation}
or equivalently the solution is the following
\[
u_+(x)=F^{-1}\{ A_{\neq}^{-1}(\xi)\tilde c(\xi')\}.
\]

Then we apply the operator $T_a$ to formula \eqref{5}
$$
(T_au_+)(t)=T_aF^{-1}\{ A_{\neq}^{-1}(\xi)\tilde c(\xi')\}
$$
and the Fourier transform
$$
(FT_au_+)(\xi)=FT_aF^{-1}\{ A_{\neq}^{-1}(\xi)\tilde c(\xi')\}.
$$

If the boundary values of our solution $u_+$ are known on $\partial C^a_+$, it means that the following function is given
$$
\int\limits_{-\infty}^{+\infty}(FT_au_+)(\xi)d\xi_m.
$$

So, if we denote
$$
\int\limits_{-\infty}^{+\infty}(FT_au_+)(\xi)d\xi_m\equiv\tilde g(\xi'),
$$
then for determining $\tilde c(\xi')$ we have the following equation
\begin{equation}\label{6}
\int\limits_{-\infty}^{+\infty}(FT_aF^{-1})\{ A_{\neq}^{-1}(\xi)\tilde c(\xi')\}d\xi_m=\tilde g(\xi'),
\end{equation}

This is a convolution equation, and if evaluating the inverse Fourier transform $\xi'\to x'$, we'll obtain the conical analogue of layer potential.

\subsection{Studying the last equation}

Now we'll try to determine the form of the operator $FT_aF^{-1}$ (see above Sec. 3). We write
\begin{equation}\label{7}
(FT_aF^{-1}\tilde u)(\xi)=(FT_au)(\xi)=\int\limits_{{\bf R}^{m-1}}e^{-ia|y'|\xi_m}e^{-iy'\cdot\xi'}\hat u(y',\xi_m)dy',
\end{equation}
where $y'=(y_1,...y_{m-1}), \hat u$ is the Fourier transform of $u$ on last variable $y_m$.

Let's denote the convolution operator with symbol $A_{\neq}^{-1}(\xi)$ by letter $a$, so that by definition
$$
(a*u)(x)=\int\limits_{{\bf R}^m}a(x-y)u(y)dy,
$$
or, for Fourier images,
$$
F(a*u)(\xi)=A_{\neq}^{-1}(\xi)\tilde u(\xi).
$$

As above let's denote $\hat a(x',\xi_m)$ the Fourier transform of convolution kernel $a(x)$ on the last variable $x_m$. The integral in \eqref{6} takes the form (according to \eqref{7})
\[
\int\limits_{{\bf R}^{m-1}}e^{-ia|y'|\xi_m}e^{-iy'\cdot\xi'}(\hat a*c)(y',\xi_m)dy',
\]

Taking into account the properties of convolution operator
and the Fourier transform we have the following representation (see Sec.3)
$$
E_a*(A_{\neq}^{-1}(\xi)\tilde c(\xi')),
$$
or, enlarged notice,
$$
\int\limits_{{\bf R}^{m-1}}E_a(\xi'-\eta',\xi_m)A_{\neq}^{-1}(\eta',\xi_m)\tilde c(\eta')d\eta'.
$$

Then the equation \eqref{6} will take the following form respectively
\begin{equation}\label{8}
\int\limits_{{\bf R}^{m-1}}K_a(\eta',\xi'-\eta')\tilde c(\eta')d\eta'=\tilde g(\xi'),
\end{equation}
where
$$
K_a(\eta',\xi')=\int\limits_{-\infty}^{+\infty}\frac{E_a(\xi',\xi_m)}{A_{\neq}(\eta',\xi_m)}d\xi_m.
$$

So, the integral equation \eqref{8} is an equation for determining $\tilde c(\xi')$. This is a conical analogue of the double layer potential.

Let's suppose we solved this equation and constructed the inverse operator $L_a$, so that $L_a\tilde g=\tilde c$. By the way we'll note the unique solvability condition for the equation \eqref{8} (i.e. existence of bounded operator $L_a$) is necessary and sufficient for unique solvability for our Dirichlet boundary value problem. Using the formula \eqref{5} we obtain
$$
\tilde u_+(\xi)=A_{\neq}^{-1}(\xi)(L_a\tilde g)(\xi'),
$$
or renaming,
\[
\tilde u_+(\xi)=A_{\neq}^{-1}(\xi)\tilde d_a(\xi').
\]

Then,
\begin{equation}\label{9}
u_+(x',x_m)=\int\limits_{{\bf R}^{m-1}}W(x'-y',x_m)d_a(y')dy',
\end{equation}
where $W(x',x_m)=F^{-1}_{\xi\to x}(A_{\neq}^{-1}(\xi))$.

The formula \eqref{9} is an analogue of Poisson integral for a half-space.

\section{Comparison with half-space case for the Laplacian}\label{sec. 7}

For the half-space $x_m>0$ we have the following (see Eskin's book \cite{E}):
$$
\tilde u_+(\xi)=\frac{\tilde c(\xi')}{\xi_m+i|\xi'|}.
$$

If we have the Dirichlet condition on the boundary, it means, that the function
$$
\tilde g(\xi')=\int\limits_{-\infty}^{+\infty}\tilde u_+(\xi)d\xi_m
$$
is given.

From formula above we have
$$
\tilde g(\xi')=\tilde c(\xi')\int\limits_{-\infty}^{+\infty}\frac{d\xi_m}{\xi_m+i|\xi'|},
$$
and we need to calculate the last integral only.

For this case we can use the residue technique and obtain, that the last integral is equal to $-\pi i$.

Thus,
$$
\tilde u_+(\xi)=-\frac{\tilde g(\xi')}{\pi i(\xi_m+i|\xi'|)}.
$$

Consequently, our solution $u_+(x)$ is the convolution (for first $(m-1)$ variables) of the given function $g(x')$ and the kernel defined by inverse Fourier transform of function $(\xi_m+i|\xi'|)^{-1}$ (up to constant).
The inverse Fourier transform on variable $\xi_m$ leads to the function $e^{-x_m|\xi'|}$, and further, the inverse Fourier transform $\xi'\to x'$ leads to Poisson kernel
$$
P(x',x_m)=\frac{c_mx_m}{(|x'|^2+x_m^2)^{m/2}},
$$
$c_m$ is certain constant defined by Euler ${\Gamma}$-function.

Thus, for the solution of the Dirichlet problem in half-space ${\bf R}^m_+$ for the Laplacian with given Dirichlet data $g(x')$ on the boundary ${\bf R}^{m-1}$
we have the following integral representation
$$
u_+(x',x_m)=\int\limits_{{\bf R}^{m-1}}P(x'-y',x_m)g(y')dy'.
$$

\section{Oblique derivative problem}\label{sec. 7}

Let' go back to formula \eqref{5}. We can write
\[
\xi_k\tilde u_+(\xi)=\xi_kA_{\neq}^{-1}(\xi)\tilde c(\xi),
\]
or equivalently according to Fourier transform properties
\[
\frac{\partial u_+}{\partial x_m}=F^{-1}\{\xi_kA_{\neq}^{-1}(\xi)\tilde c(\xi)\},
\]
for arbitrary fixed $k=1,2,...,m$.

Further, we apply the operator $T_a$ and work as above. Our considerations will be the same, and in all places instead of $A_{\neq}^{-1}(\xi)$ will stand $\xi_kA_{\neq}^{-1}(\xi)$.

I call this situation the oblique derivative problem, because $\frac{\partial}{\partial x_k}$ related to conical surface is not normal derivative exactly.

{\bf Remark 2.}
{\it  Some words on Neumann problem. If we try to give normal derivative of our solution on conical surface different from origin, then we have the boundary value problem with variable coefficients because the
boundary condition varies from one point to another one on conical surface. We need additional localization for such points to reduce it to the case of constant coefficients and consider corresponding model problem
in ${\bf R}^m_+$. Roughly speaking, I would like to say, that the solution looks locally different in dependence on the type of boundary point. In other words, local principle permits to work with symbols and boundary
conditions non-depending on space variable.}

\section{Conclusions}

It seems to solve explicitly the simplest boundary value problems in domains with conical point we need to use another potentials different from classical
simple and double layer potentials. I will try to show this fact for the Laplacian with Dirichlet condition on conical surface in my forthcoming paper
by direct calculations.

{\bf Acknowledgements.}
Many thanks to DAAD and Herr Prof. Dr. Volker Bach for their support.

\vspace{2cm}
Chair of Pure Mathematics, Lipetsk State Technical University,\\
Moskovskaya 30, Lipetsk 398600, Russia.\\
E-mail: vbv57@inbox.ru

\end{document}